\documentclass[sn-mathphys]{sn-jnl}

\jyear{2021}%

\theoremstyle{thmstyleone}%
\theoremstyle{thmstyletwo}%

\usepackage{changepage}
\usepackage{fullwidth}

\theoremstyle{thmstylethree}%

\begin{document}

\title[    ]{NEURAL NETWORK ALGORITHM FOR INTERCEPTING TARGETS MOVING ALONG KNOWN TRAJECTORIES BY A DUBINS' CAR%
\thanks{The work was supported by a grant from the ICS RAS Youth Scientific School. <<Methods of optimization and motion planning of controlled objects>>.
The work of A.A. Galyaev and I.A. Nasonov was partially supported by the Russian Scientific Foundation (project №~23-19-00134).}}


\author*[1]{\fnm{Ivan} \sur{Nasonov}}\email{nasonov.ia18@physics.msu.ru}
\equalcont{These authors contributed equally to this work.}

\author[1]{\fnm{Andrey} \sur{Galyaev}}\email{galaev@ipu.ru}
\equalcont{These authors contributed equally to this work.}

\author[1]{\fnm{Andrey} \sur{Medvedev}}\email{medvedev.ai18@physics.msu.ru}
\equalcont{These authors contributed equally to this work.}

\affil*[1]{\orgdiv{Laboratory 38}, \orgname{Institute of Control Sciences of RAS},\\ \orgaddress{\street{} \city{Moscow}, \postcode{}\state{}\country{Russia}}}


\abstract{The task of intercepting a target moving along a rectilinear or circular trajectory by
a Dubins' car is formulated as a time-optimal control problem with an arbitrary
direction of the car's velocity at the interception moment. To solve this problem
and to synthesize interception trajectories, neural network methods of unsupervised
learning based on the Deep Deterministic Policy Gradient algorithm are used. The
analysis of the obtained control laws and interception trajectories in comparison
with the analytical solutions of the interception problem is performed. The
mathematical modeling for the parameters of the target movement that the neural
network had not seen before during training is carried out. Model experiments are
conducted to test the stability of the neural solution. The effectiveness of using
neural network methods for the synthesis of interception trajectories for given
classes of target movements is shown.
}

\keywords{Interception task, Dubins' car, DDPG algorithm, neural network synthesis of trajectories.} 



\maketitle

\section{Introduction}

The task of intercepting mobile targets moving along known trajectories has been of interest to researchers since the mid-50s of the last century \cite{isaacs}. One of the basic models for describing the dynamics of an intercepting object is the Dubins' car model.

The first works on finding a line with a limited curvature and a minimum length connecting two given points belong to A.A. Markov. His first task in \cite{markov} was devoted to finding a curve connecting two points on a plane with minimal length and bounded curvature with a fixed exit direction from the first point. Such a task has found application in solving the problems of laying railways. In 1957, L. Dubins published a similar work \cite{dubins} on finding a line of the shortest length with a limited radius of curvature connecting two points on a plane with a given direction of exit from the first point and a given direction of entry into the second. The results proved to be useful in the study of objects with a limited turning radius and a constant speed of movement.

In \cite{control} the non-game problem of the fastest interception of a moving target by a Dubins' car is considered. It was assumed that the target was moving along an arbitrary and previously known continuous trajectory. To find the solution, the algebraic criterion of the optimality of the interception along the geodesic line and the optimal value of the interception time criterion were found.

In early studies of \cite{trajectories}, sufficient conditions were established that the optimal trajectory is curves. These conditions impose restrictions on the ratio of the minimum radius of curvature of the trajectory of the car and the distance between the target and the car at the initial moment of time. In \cite{berdishev}, control has been synthesized to intercept a target along a geodesic line drawn from the beginning of the movement of the car to the intercept point, and it is assumed that the target is moving in a straight line with a constant speed.

The practical applications of the tasks of interception by the Dubins' car are quite extensive: the construction of optimal trajectories of unmanned aerial vehicles that monitor several ground targets \cite{planes}, the development of algorithms that solve the traveling salesman problem \cite{comiv}, the construction of bypass trajectories when moving with obstacles \cite{path}. Also, the Dubins' car model is used in the pursuit-evasion differential game. Such a game involves the presence of two agents: the pursuer must catch the target, and the escapee must evade the pursuer. An analytical solution to the problem of finding the optimal interception time and synthesis of the optimal trajectory for such a game was obtained in \cite{control}. The problem of synthesis of intercept trajectories for objects moving along a circular trajectory was considered in \cite{manyam}.

The solution of the problems of interception by the Dubins' car can also be obtained with the help of computers. Recently, neural network reinforcement learning methods have been actively used for such tasks, which represent machine learning technology without models and are used in cases when there is little or no data for training a neural network at all. Unlike learning with a teacher \cite{sl}, who needs a set of marked-up data, reinforcement learning is based on the interaction of the agent with the environment \cite{rl}. This method is most effective for finding a solution to the problem of pursuit-evasion.
	
The Actor-Critic method is used in many relevant studies. For example, in \cite{Perot}, Actor-Critic was used with Convolutional Neural Network (CNN) and Long Short-Term Memory (LSTM) as a state encoder for racing games. In \cite{maolinwang}, a fuzzy deterministic policy gradient algorithm was used to obtain a specific physical meaning when teaching politics in the pursuit-evasion game. In \cite{Lillicrap}, the Deep Deterministic Policy Gradient (DDPG) method for interacting with a continuous action space was introduced for the first time. It is this algorithm that will be used in this work for neural network synthesis of the trajectory of interception by the Dubins' car of a target moving at a constant speed along rectilinear and circular trajectories. Thanks to DDPG, it was possible for the first time to obtain a suboptimal trajectory based on a neural network \mbox{solution}.

The relevance of the work is due to both the demand in practice for interception algorithms for one and many moving targets, and the possibility of obtaining some new theoretical results related to the synthesis of interception trajectories. Of particular interest is the so-called traveling salesman problem with mobile goals --- Moving Target Traveling Salesman Problem (MTTSP) \cite{MTTSP}. In this case, the points that need to be bypassed are moving at a given speed. An example of such a scenario is the interception of several evading (or attacking) targets, which is very important for dual-use applications. Obviously, finding the best route to intercept several mobile targets is a particularly difficult task due to the constant change in the position of targets, which significantly increases the computational costs of finding optimal solutions. It is known that a heuristic approach has been proposed in the literature to solve MTTSP.

The authors propose a synthesis of the interception trajectory based on a neural network solution, since analytical results and optimal trajectories for groups of targets are practically absent or unknown. The authors plan to scale this method for similar tasks.

The structure of the work includes 6 sections. Section 2 offers a mathematical formulation of the problem adapted for further application. Section 3 is devoted to the description of the DDPG algorithm, also ready for use in this formulation. Section 4 describes the structure of the neural network, and section 5 contains the simulation results. In conclusion, the direction of further research is presented.

\section{Formulation of the neural network interception problem}\label{s:sec2}

On the plane, the problem of the fastest $\delta$ is considered-the interception by the Dubins' car (pursuer) of a moving object (target) moving along two given trajectories at a constant speed. As in \cite{control}, the dynamics for the pursuer was selected as
    \begin{gather}\label{Eq1}
    \left\{
    \begin{array}{l}
   \displaystyle \dot{x_P}=\cos{\varphi},
    \\
    \displaystyle\dot{y_P}=\sin{\varphi}, \\
    \dot{\varphi}=u,~ |u(t)|\leqslant1.
    \end{array}
    \right.
    \end{gather}

    \begin{figure}[H]
		
		\centering
		
		\includegraphics[width=0.5\linewidth]{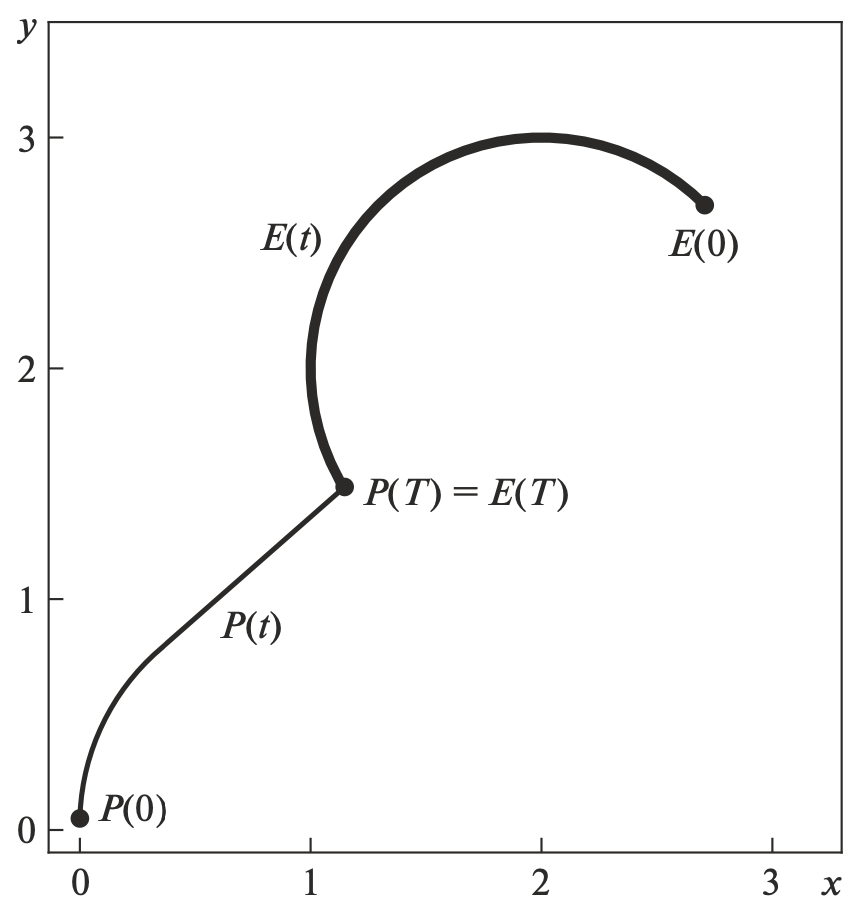}
		
		\caption{Mutual location of objects.}
		\label{fig.problemstatement}

    \end{figure}

Here $x_P(t)$ and $y_P(t)$ are the coordinates of the Dubins' car on the Cartesian plane, $\varphi(t)$ is the angle between the direction of the pursuer's speed and the abscissa axis, and $u(t)$ is a time-dependent control that shown in Fig.\;\ref{fig.problemstatement}. The coordinates and angle of the car are denoted by the vector function $P(t) = (x_P(t),y_P(t), \varphi(t))$.

The initial conditions of the system (\ref{Eq1}) are fixed:
    \begin{gather}\label{Eq2}
    x_P(0) = 0, \quad y_P(0) = 0, \quad \varphi(0) = \frac{\pi}{2}.
     \end{gather}

Continuous vector function $E(t)=(x_E(t),y_E(t))$ defines the trajectory of the target on the Cartesian plane.

The terminal condition of $\delta$-interception for a neural network solution has the following form:
    \begin{gather}\label{Eq3}
    (x_P(T) - x_E(T))^2 + (y_P(T) - y_E(T))^2 \leqslant \delta^2,
     \end{gather}
where $T \in \mathbb{R}^{+}_{0}$~--- the time of movement from the starting point to the interception point, and ~$\delta$~--- the specified interception radius~--- the maximum allowable distance between the pursuer and the target at which the interception it can be considered perfect. This parameter is introduced to define the concept of interception specifically for a neural network solution.

Let's set the task of intercepting the target in minimal time as an optimal control problem in the class of piecewise constant functions:
    \begin{gather}\label{Eq4}
    \displaystyle J[u] \stackrel{def}{=} \int\limits_0^T dt \rightarrow \underset{u}{\min}.
    \end{gather}

Let's start describing the dynamics of the goal. According to the condition of the task, the target moves at a constant speed in a straight line or in a circle. Then the parametrized coordinate equations will have the following form:
    \begin{gather}\label{Eq5}
    \begin{cases}
    x_{E}(t)= R\cos(\omega t+\phi)+x_{0},
    \\
    y_{E}(t)=R\sin(\omega t+\phi)+y_{0};
    \end{cases}
\\[.6em]
\label{Eq6}
    \begin{cases}
    x_E(t)=v_x t + x_0,
    \\
    y_E(t)=v_y t + y_0,
    \end{cases}
   \end{gather}
where $x_0$ and $y_0$ are the initial conditions of the target coordinates and are chosen arbitrarily.

To take into account the relative position of the pursuer and the target, we introduce a formula for finding the angle between the abscissa axis and the straight line connecting the coordinate points of the target and the pursuer. Let $(x_P, y_P)$ and $(x_E, y_E)$~--- the coordinates of the pursuer and the target, respectively, at some point in time $t$. Then the desired value of the angle is found by the formula
    \begin{gather*}
    \psi = \arctan{ \left( \frac{y_E - y_P}{x_E - x_P} \right)}.
    \end{gather*}
    We will also introduce a formula for calculating the distance $L$ between agents:
    \begin{gather*}
    L=\sqrt{(x_P-x_E)^2 + (y_P-y_E)^2}.\end{gather*}

Next, to simplify the study of the problem, we will make the transition to the new coordinates. To do this, you need to be able to compare the current state of agents ${S=(x_P, y_P, \varphi, x_E, y_E)}$ and the state predicted by the neural network ${S'=(x'_P, y'_P, \varphi', x'_E, y'_E)}$.

We get the values for the functions of the angles $\psi$ and $\psi'$ from the states $S$ and $S'$, respectively, and also calculate the distance $L'$ when the agents are in the state~$S'$. We introduce the angle between the direction of the speed of the pursuer and the line connecting the coordinate points of the agents:
    \begin{gather*}
    \Theta=\varphi' - \psi'.\end{gather*}
    Let's introduce the rotation speed as a quotient of the difference ${\psi' -\psi}$ and the time interval~$\Delta t$ during which the transition from the state $S$ to the state $S'$ occurred:
    \begin{gather*}
    \omega = \frac{\psi' - \psi}{\Delta t}.\end{gather*}

    The totality of $(L', \omega, \Theta)$ and there are the desired coordinates in which we will build a neural network solution.
    At the initial moment of time, when the result of the neural network has not yet been received, the coordinates are $(L'(0), \omega(0), \Theta(0))$ are calculated as follows:
    \begin{gather}\label{Eq7}
    \begin{cases}
    L'(0)=L(S),
    \\
    \omega(0) = 0,
    \\
    \Theta(0) = \varphi(0) - \psi(0),
    \end{cases}
    \end{gather}
    where
   $ \psi(0) = \arctan{ \left( \frac{y_E(0) - y_P(0)}{x_E(0) - x_P(0)}
    \right)}$. 

\section{Algorithm Deep Deterministic Policy Gradient} \label{s:methods}
DDPG~--- is an Actor-Critic algorithm based on a deterministic policy gradient. The DPG (Deterministic Policy Gradient) algorithm consists of a parameterized function Actor $\mu\left(s\mid \theta^{\mu}\right)$, which sets control at the current time by deterministic matching of states with a specific action. The function Critic $Q(s,a)$ is updated using the Bellman equation in the same way as with $Q$ training. The Actor is updated by applying a chain rule to the expected reward from the initial distribution of $J$ in relation to the parameters of the Actor:

\begin{equation}
\begin{array}{ccc}   

\displaystyle \nabla \theta^{\mu} \mathrm{J} \approx \mathbb{E}_{s_{\mathrm{t}} \sim \rho^{\beta}}\left[\left.\nabla_{a} Q\left(s, a \mid \theta^{Q}\right)\right \vert_{s=s_{\mathrm{t}}, a=\mu\left(s_{\mathrm{t}} \mid \theta^{\mu}\right)}\right] \\
\displaystyle =\mathbb{E}_{\mathrm{s_t} \sim \rho^{\beta}}\left[\left.\left.\nabla_{\theta^{\mu}} Q\left(s, a \mid \theta^{Q}\right)\right\vert_{s=s_{\mathrm{t}}, a=\mu\left(s_{\mathrm{t}}\right)} \nabla_{\theta_{\mu}} \mu\left(s \mid \theta^{\mu}\right)\right\vert_{s=s_{\mathrm{t}}}\right].

\end{array}
\end{equation}


DDPG combines the advantages of its predecessors, which makes it more stable and effective in training. Since different trajectories can be very different from each other, DDPG uses the idea of DQN \cite{dqn}, called a playback buffer. The playback buffer~--- is a finite-size buffer into which media data is stored at any given time. It is necessary to achieve a uniform distribution of the transition sample and discrete control of neural network training. Actor and Critic are updated by evenly sampling the mini-batch from the playback buffer. Another addition to DDPG was the concept of updating program targets instead of directly copying weights to the target network. Network being updated $Q\left(s, a\mid\theta^{Q}\right)$ is also used to calculate the target value, so updating $Q$ is subject to divergence.This is possible if you make a copy of the Actor and Critic networks, $Q^{\prime}\left(s, a\mid \theta^{Q^{\prime}}\right)$ and $\mu^{\prime}\left(s, a\mid \theta^{\mu^{\prime}}\right)$. The weights of these networks are as follows: ${\theta^{\prime} \leftarrow \tau\theta+(1-\tau) \theta^{\prime}}$
with $\tau\ll 1$. The research problem is solved by adding the noise received from the noise process $N$ to the control of the actor. In this study, the Ornstein-Uhlenbeck process is selected
\cite{Uhlenbeck}.

\begin{figure}[t]
\centering{\includegraphics[width=130mm]{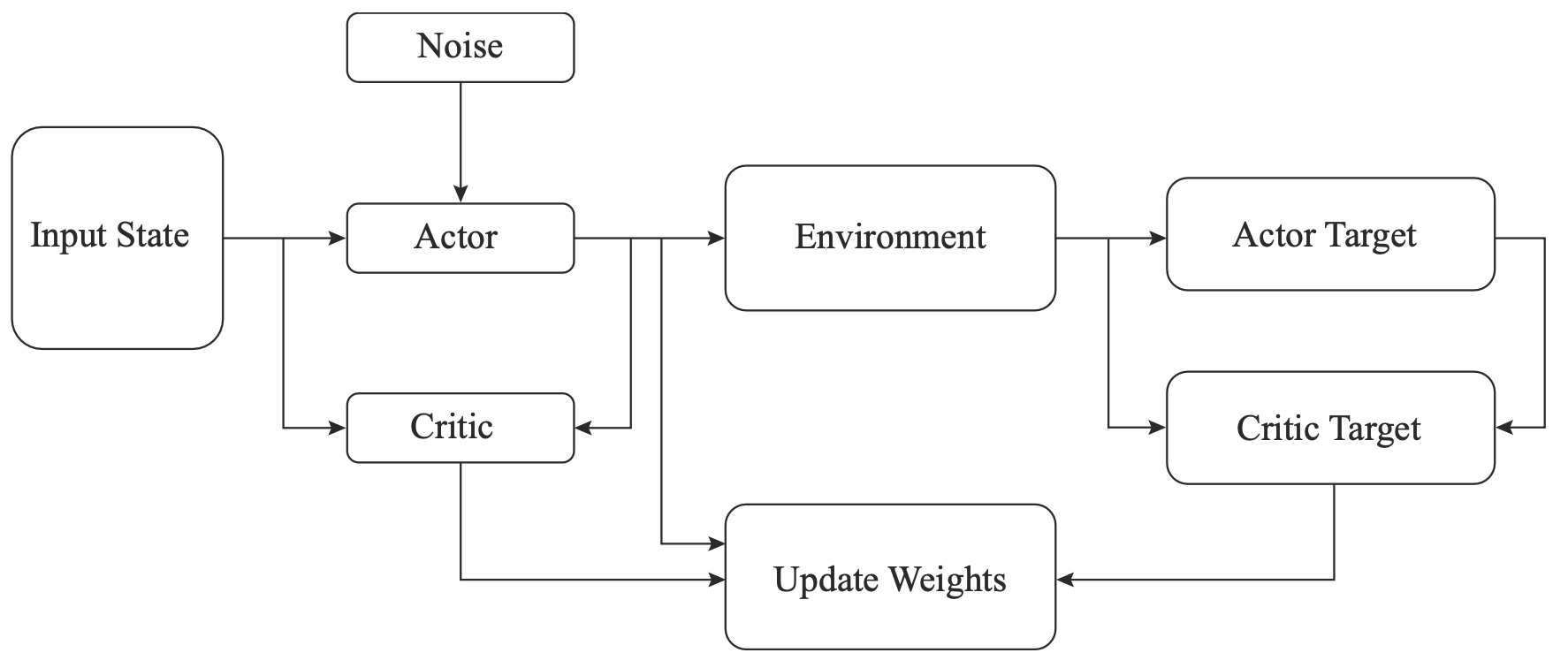}}
		\caption{The general structure of the Deep Deterministic Policy Gradient algorithm.}
		\label{fig: ddpg}
\end{figure}

The general structure of DDPG is shown in Fig.\;\ref{fig: ddpg}. Since the task requires that the controls are enclosed in a numerical interval, it is necessary to introduce restrictions. To do this, the program used the $clip()$ function, which limits the range of action values in the range $[-1;1]$.

    \begin{algorithm}[t]
    \caption{Deep Deterministic Policy Gradient}\label{alg:DDPG}
    \hspace*{\algorithmicindent} \textbf{Input Data:} discount coefficient $\gamma$, number of episodes $M$, number of training steps $T$ in each episode, batch size $N$, training coefficients of neural networks Actor and Critic $r_a$ and $r_c$, respectively.\\
    \hspace*{\algorithmicindent} \textbf{Output Data:} Control $u=\mu(s\vert\theta^\mu)$
    \begin{algorithmic}[1]

    

    \State Arbitrary initialization of the networks Actor $\mu(s\vert\theta^\mu)$ and Critic $Q(s, a\vert\theta^Q)$

    \State Initialization of target networks $Q'$ and $\mu'$ with weight parameters $\theta^Q=\theta^{Q'}$ and~$\theta^\mu=\theta^{\mu'}$

    \State Initializing the $R$ buffer
    \For{ $episode \gets 1$ to $M$}

    \State Initialization of a random action $a_t=\mu(s_t\vert\theta^\mu)+\eta_t$ according to the current control and research noise

    \State Getting the initial state of the $s_1$ environment
    \For{ $t \gets 1$ to $T$}

    \State\;\,\,Performing the action $a_t$, acquiring the reward $r_t$ and obtaining a new state of the environment $s_{t+1}$

    \State\;\,\,Saving the transition $(s_t, a_t, r_t, s_{t+1})$ in the buffer $R$

    \State Random sampling of $N$ transitions $(s_i, a_i, r_i, s_{i+1})$ from $R$

    \State Getting $y_i=r_i+\gamma Q'(s_{i+1}, \mu'(s_{i+1}\vert\theta^{\mu'})\vert\theta^{Q'})$

    \State Updating the weights of the Critic network by minimizing the loss function 
    $$\widehat{L}=\frac{1}{N}\sum_{i}(y_i-Q(s_i, a_i\vert\theta^Q))^2$$

    \State Updating an Actor Policy with an Effective Policy Gradient: 
    $$\nabla_{\theta^{\mu}}J\approx \frac{1}{N} \sum_{i} \nabla_{a} Q(s, a\vert\theta^Q) \vert_{s=s_i, a=\mu(s_i)} \nabla_{\theta^\mu} \mu(s\vert\theta^\mu)\vert_{s_i}$$

    \State Updating target networks \begin{gather*}
    \theta^{Q'}=\tau \theta^{Q} + (1-\tau)\theta^{Q},\quad \theta^{\mu'}=\tau \theta^{\mu} + (1-\tau)\theta^{\mu}\end{gather*}
    
    \EndFor
    \EndFor
   \end{algorithmic}
\end{algorithm}

Table 1 shows the differences between the Actor, Critic networks and their target networks. It contains input and output values, as well as formulas for calculating these values.

A detailed description of the DDPG method is given in the algorithm
\ref{alg:DDPG}.

\begin{table}[t]
    \label{table:1}
    \resizebox{1.1\textwidth}{!}{\begin{tabular}[H]{|p{18mm}|p{35mm}|p{35mm}|p{35mm}|}
    \hline
    \multicolumn{1}{|c|}{\rule{0pt}{4mm}Network}
     & \multicolumn{1}{c|}{Formula} & \multicolumn{1}{c|}{Input Data} & \multicolumn{1}{c|}{Output data} \\ [0.25ex]
    \hline
    \rule{0pt}{4mm} Critic target & $Q'(s_{t+1}, \mu ' (s_{t+1} \vert \theta^{\mu'} )\vert\theta^{Q'})$ & the next state of the environment; the output of the target network Actor & \mbox{value $Q'$,} which is used to calculate~$y_i$
    \\
    \rule{0pt}{4mm} Critic & $ Q(s_t, a \vert \theta^{Q} ) $ & current state of the environment; current action & the $Q$ value that is needed to calculate the loss and update the Actor network
    \\
    \rule{0pt}{4mm}Actor target & $\mu'(s_{t+1} \vert \theta^{\mu'})$ & the next state of the environment & the action $\mu'$ used as the input value of the target network Critic      \\
    \rule{0pt}{4mm}Actor & $\mu(s_t \vert \theta^u)$ & current state of the environment & the $\mu$ action that is used to update the Actor network
    \\
    \hline
    \end{tabular}}

\end{table}

\section{Neural network}
\subsection{Network Architecture}
\begin{figure}[t]
\centering{\includegraphics[width=65mm]{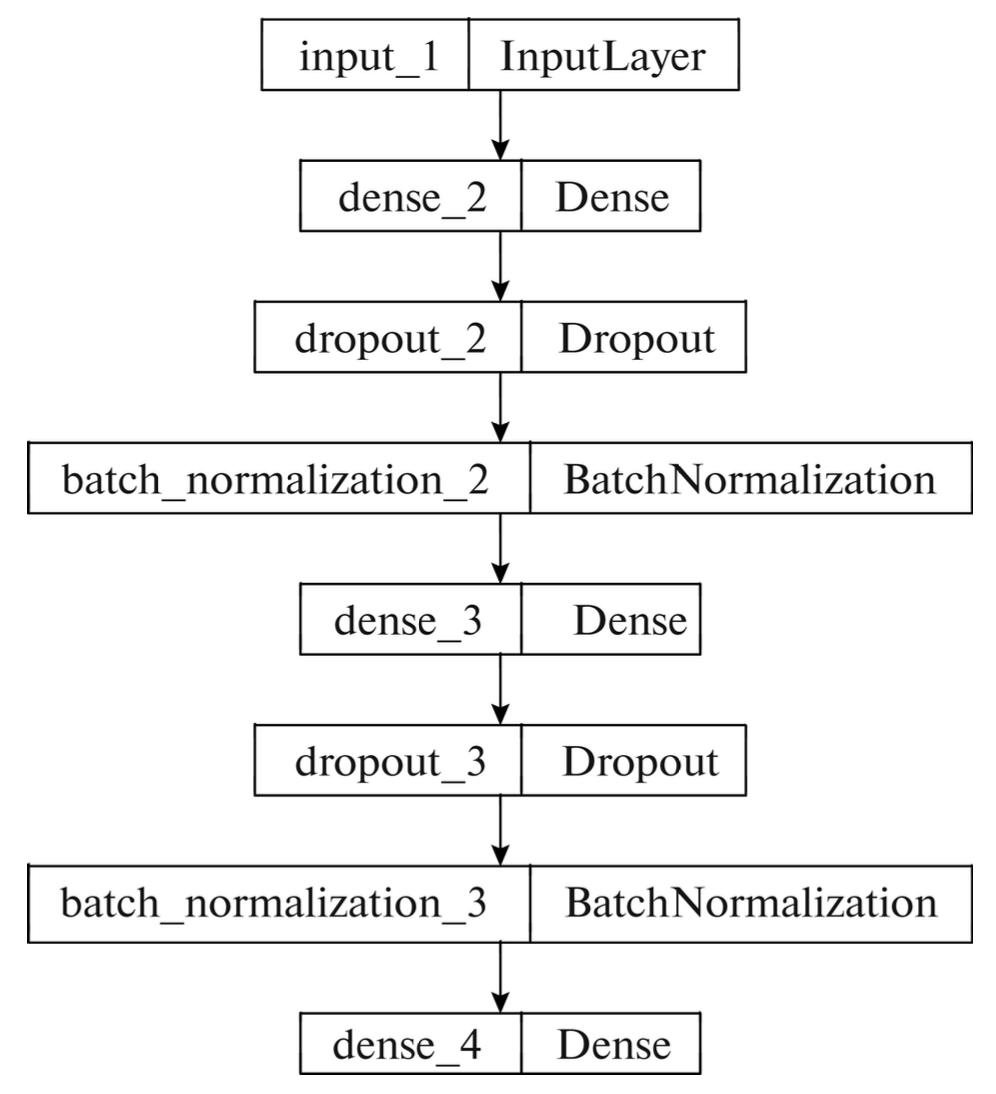}}
		\caption{Actor neural network architecture.}
		\label{fig:actor_model}
    \end{figure}

\begin{figure}[t]
\centering{\includegraphics[width=100mm]{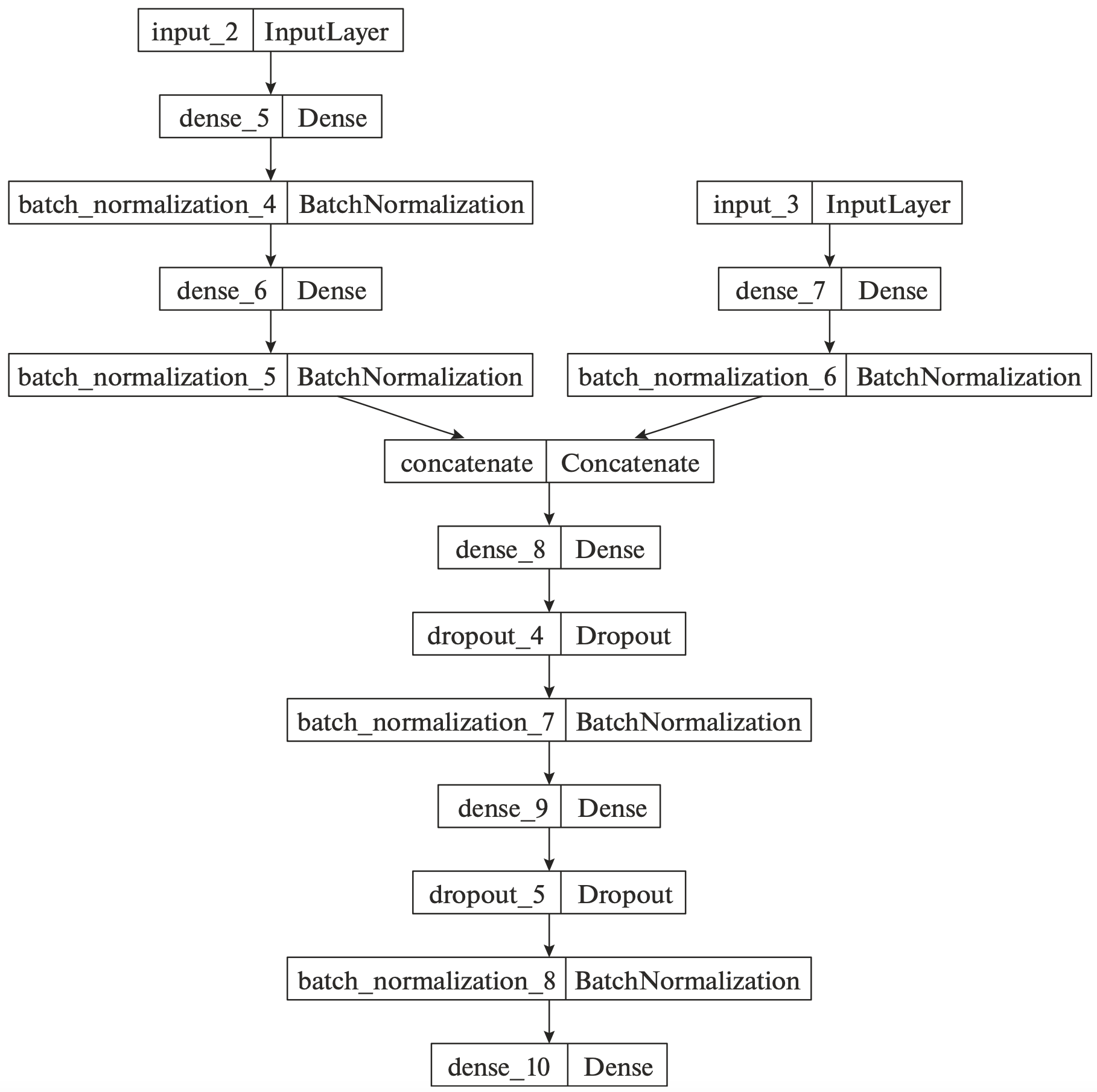}}
		\caption{The architecture of the neural network Critic.}
		\label{fig:critic_model}
    \end{figure}

To implement the Deep Deterministic Policy Gradient algorithm, two neural networks were written for each method: Critic and Actor. Their architectures are depicted in Fig. \ref{fig:critic_model} and \ref{fig:actor_model}.

The Actor network has four fully connected hidden layers with 256 neurons, with ~$SELU$ activation function. Since the possible actions are in the range $[-1,1]$, it is convenient to take the activation function for the output layer as $tanh$. The Critic network has five fully connected hidden layers with 16, 32, 32 and two layers with 512 neurons, with an activation function $SELU$.

The Critic and Actor networks are made up of fully connected $Dense$ layers, for the output values of which the normalization operation and the $Dropout$ \cite{dropout} method are used, which is effective in combating the problem of retraining neural networks. To calculate the output of the Actor network from the last layer, the hyperbolic tangent activation function is selected.

The Critic network has a complex structure because it takes two input values: the state of the environment and the actions of the pursuer. Next, the layers are connected using the $Concatenate$ method and the values pass through the fully connected layers of the network to the output, which is a layer of unit dimension.

\subsection{Hyperparameters}

    The $SELU$ \cite{selu} function was chosen as the activation function in the hidden layers of the Critic and Actor neural networks, which is given by the following equation:
\begin{gather*}
    SELU(x) = \lambda
    \begin{cases}
    x, & x>0
    \\
    \alpha e^x - \alpha, & x\leqslant0,
    \end{cases}
    \end{gather*}
    where $\lambda\approx1{,}0507$, and $\alpha\approx1{,}6732$.

\begin{figure}[t]
\centering{\includegraphics[width=65mm]{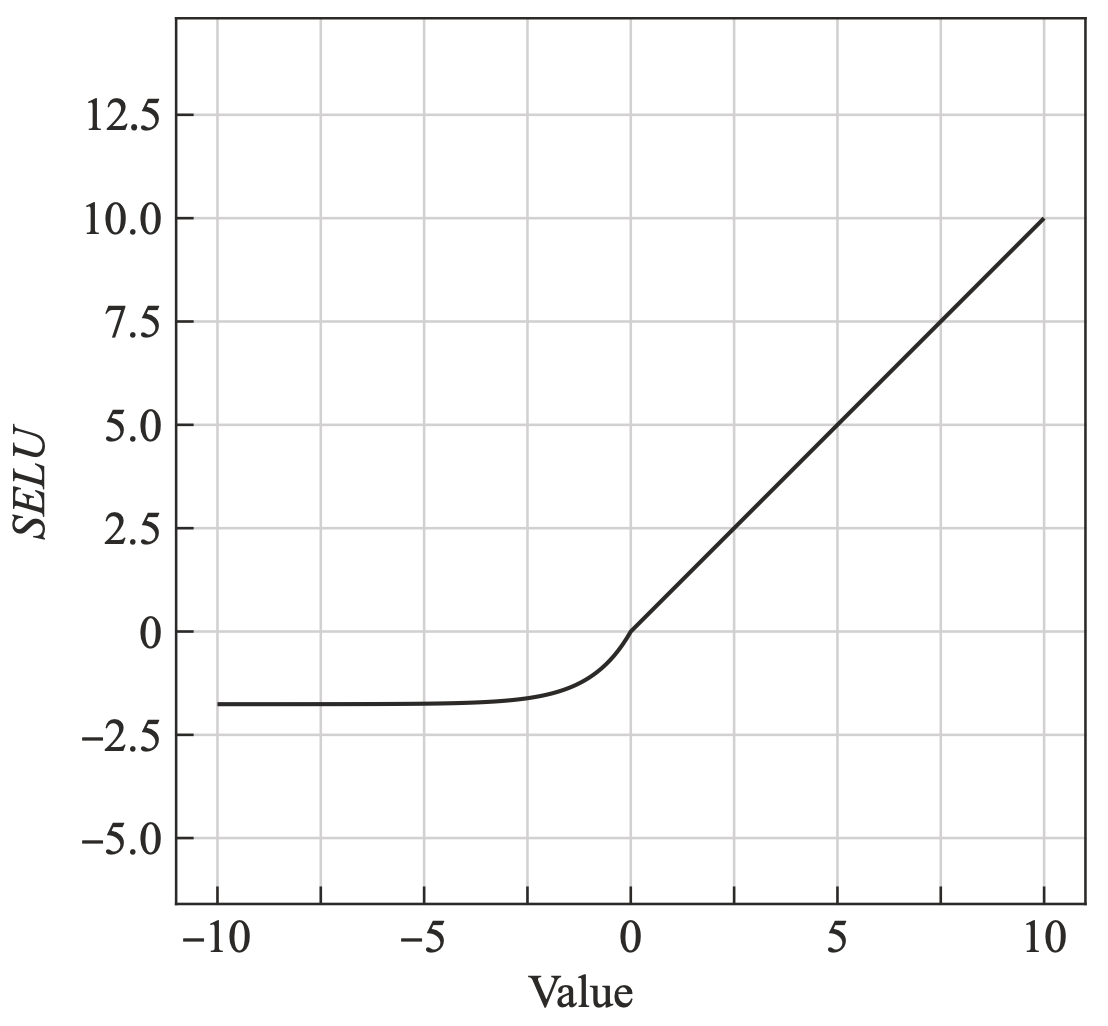}}
		\caption{Graph of the activation function $SELU$.}
		\label{fig:selu}
\end{figure}

    The graph of the $SELU$ function is shown in Fig. \ref{fig:selu}.

    The $SELU$ function has the property of self-normalizing input data when using the $LeCun$ initialization method, which initializes network parameters as a normal distribution. Therefore, the output values of this function have a zero mean and a single standard deviation.

\begin{figure}[t]
\centering{\includegraphics[width=65mm]{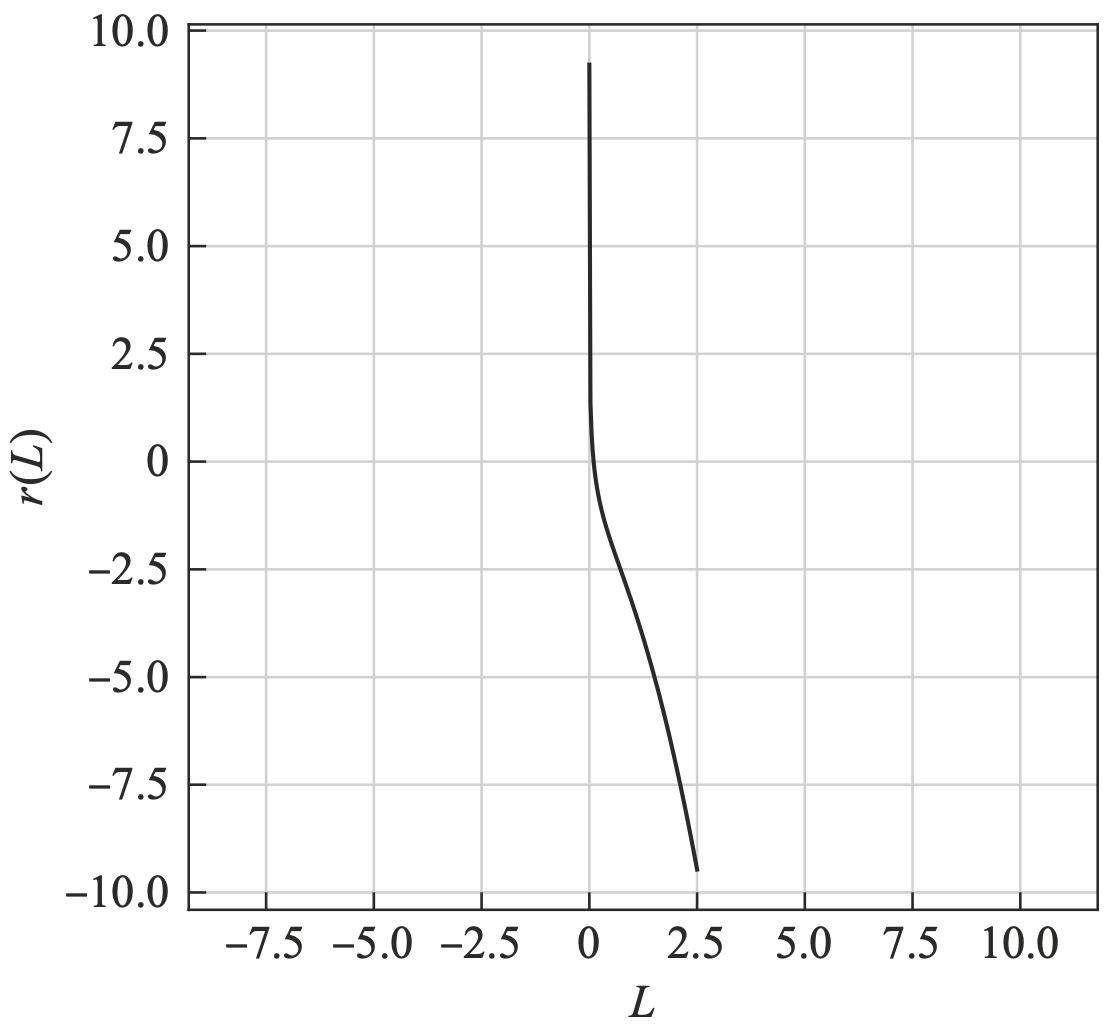}}
		\caption{A graph of the dependence of remuneration on the distance between agents.}
		\label{fig:rew_view}
\end{figure}

    In the form of a reward function for the pursuer, the following expression was chosen, depending only on the distance $L$ between the agents:
    \begin{gather} \label{eq1}
    r(L) = -\lg{(10 L)} - L^2.
    \end{gather}
    The graph of this function is shown in Fig. \ref{fig:rew_view}. On it you can see that the value of $r$ grows rapidly with a decrease in $L$, and when the distance takes a zero value, the agent receives the maximum reward.
    
The values of hyperparameters of neural networks are given in Table \ref{table:2}. The parameters $\gamma$, $\tau$, episode size and time interval were selected as a result of the analysis in accordance with \cite{maolinwang}. However, the values of the mini-batch size, buffer volume R, step size and training coefficients of Actor-Critic networks were selected empirically~--- the network synthesized the trajectories of intercepting the movement of the target, and then their analysis was carried out for compliance with the physical task. For example, if the average reward schedule did not increase during 100-200 training episodes, and the values of the error functions of the Actor-Critic neural networks did not decrease over the same period, then the values of the training coefficients of the networks decreased, and the size of the mini-batch increased.

\begin{table}[t]
    \caption{Values of neural network parameters\hspace*{50mm}}
    \label{table:2}

    \begin{tabular}[H]{|p{35mm}|c|p{75mm}|}
    \hline
    \multicolumn{1}{|c|}{\rule{0pt}{4mm}Parameter}
     & \multicolumn{1}{c|}{Value} & \multicolumn{1}{c|}{Description} \\ [0.25ex]
    \hline
    \rule{0pt}{4mm}$\gamma$    & 0,98  & The discount factor used in the Bellman equation
    \\
    \rule{0pt}{4mm}$\tau$    & 0,01   & Coefficient of soft updating of target networks
    \\
    \rule{0pt}{4mm}Size mini-batch & 64  & Number of samples to update the weights
    \\
    \rule{0pt}{4mm}$R$ Buffer Size   & 10\,000  & The amount of data from which examples are selected for updating      \\
    \rule{0pt}{4mm}Episode Size & 1000 & Number of episodes used for training
    \\
    \rule{0pt}{4mm}Step Size & 400 & The number of training steps in each episode
    \\
    \rule{0pt}{4mm}Time interval & 0,1 & Time of each step of training
    \\
    \rule{0pt}{4mm}The learning coefficient of the Actor network & 5e-5 & The learning factor used to update the Actor network
    \\
    \rule{0pt}{4mm}The learning coefficient of the Critic network & 1e-4 & The learning factor used to update the Critic network
    \\
    \hline
    \end{tabular}

    \end{table}

\section{Simulation results}

\subsection{Neural network learning process}

The simulation was performed using Python and the TensorFlow framework. The initial parameters of the movement of the target and the pursuer during neural network training are given in Table \ref{table:3}.

    \begin{figure}[t]
\centering{\includegraphics[width=65mm]{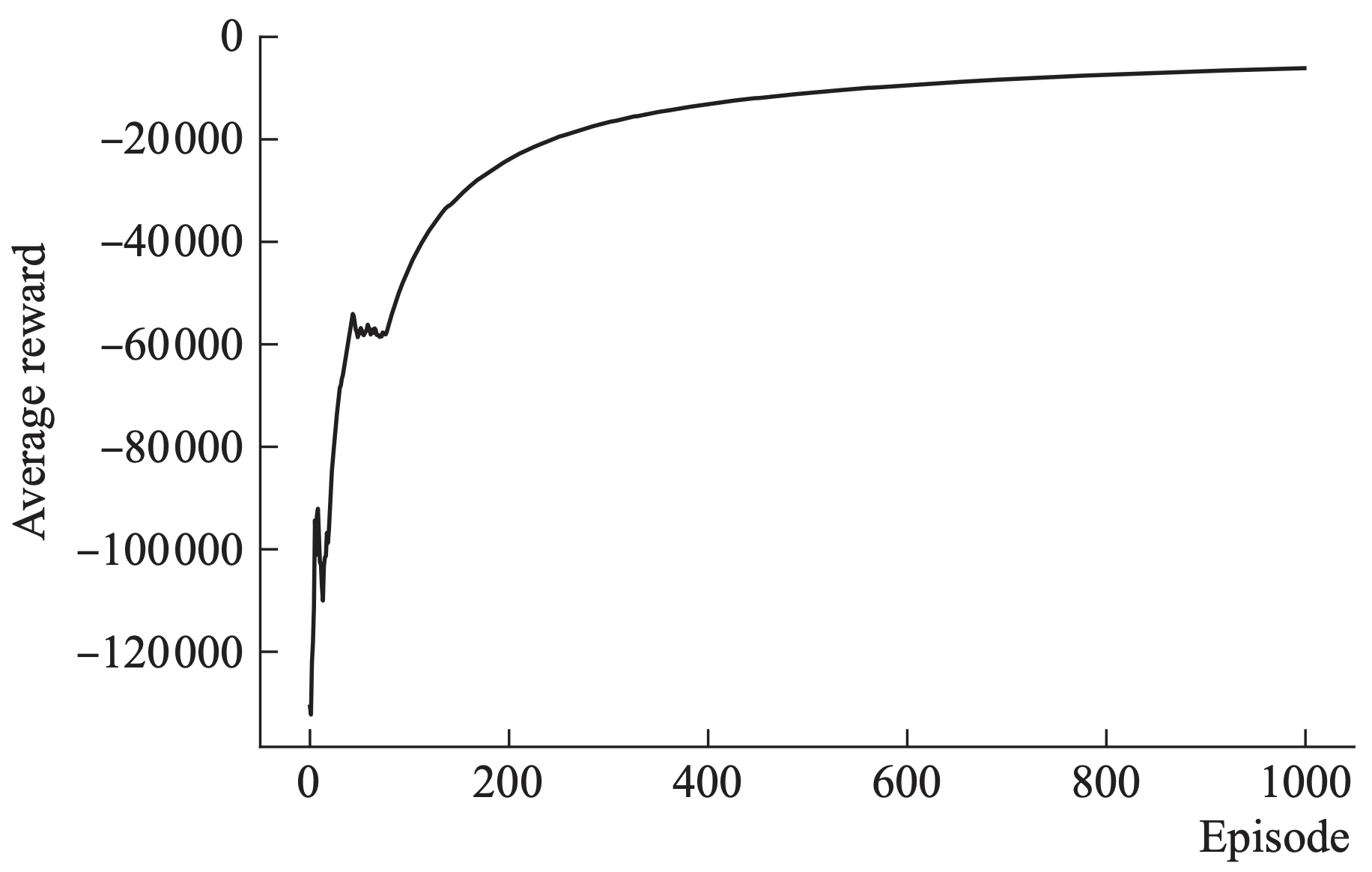}}
		\caption{The dependence of the average reward on the episode number.}
		\label{fig:avg_reward}
\end{figure}

\begin{table}[b]
    \caption{Initial parameters of target and pursuer movement during network training\hspace*{100mm}}
    \label{table:3}

\begin{center}
    \begin{tabular}[H]{|p{85mm}|p{38mm}|}
    \hline
    \multicolumn{1}{|c}{\rule{0pt}{4mm}Parameter}
     & \multicolumn{1}{|c|}{Value} \\ [0.25ex]
    \hline
    \rule{0pt}{4mm}{The initial coordinate of the target movement $x_E(0)$} &
    {An arbitrary value in the interval} $(-3;3)$
    \\
    \rule{0pt}{4mm}{The initial coordinate of the target movement $y_E(0)$} &
    {An arbitrary value in the interval} $(-3;3)$
    \\
    \rule{0pt}{4mm}{\mbox{Initial coordinates of the pursuer's movement}
    $(x_P(0);y_P(0))$} &
    $(0;0)$
    \\
    \rule{0pt}{4mm}{Initial orientation of the pursuer} $\varphi(0)$
     & $\displaystyle\frac{\pi}{2}$
    \\
    \rule{0pt}{4mm}{Constant speed of the pursuer $v$} & 1
    \\
    \rule{0pt}{4mm}{Intercept radius $\delta$} & 0,2
    \\
    \hline
    \end{tabular}
\end{center}

\vspace*{2mm}
    \caption{Characteristics of the equipment where the network was trained}
    \label{table:cpu}

    \begin{center}

    \begin{tabular}[H]{|p{7cm}|c|}
    \hline
    \multicolumn{1}{|c|}{\rule{0pt}{4mm}Parameter} & {Value} \\ [0.25ex]
    \hline
    \rule{0pt}{4mm}{Processor} & Intel(R) Core(TM) i7-8565U
    \\
     {Lithography} & 14 nm
    \\
     {Number of cores} & 4
    \\
     {Number of threads} & 8
    \\
     {Processor base clock frequency} & 1,80 hHz 
    \\
     {Cache memory} & 8 Mb
    \\
     {Computer RAM} & 16 Gb
    \\
    \hline
    \end{tabular}  \end{center}
    \end{table}

The initial coordinates of the target movement are randomly selected using the numpy.random.uniform() function in the range $(-3;3)$ so that the network trains on different examples and works effectively after the training process. The target speeds have always had a constant value throughout the learning process $v_x = 0{,}5$ and $v_y = 0{,}5$.

Neural network training was carried out on a process with the characteristics specified in Table. \ref{table:cpu}. Due to the complexity of the neural network model, the learning process lasted about four hours.

In Fig. \ref{fig:avg_reward} shows a graph of the average remuneration for the entire training period. During the training of the model, there is a sharp increase in the value of the agent's reward in the first 100--150 episodes. Filling of the playback buffer $R$ corresponds to this process. Next, the training examples are randomly taken from $R$, the network training process takes place and the resulting tuple of states replaces the old data sample in $R$. At this stage, there is a slow increase in the average remuneration, see Fig. \ref{fig:avg_reward}.

    \begin{figure}[t]
\centering{\includegraphics[width=65mm]{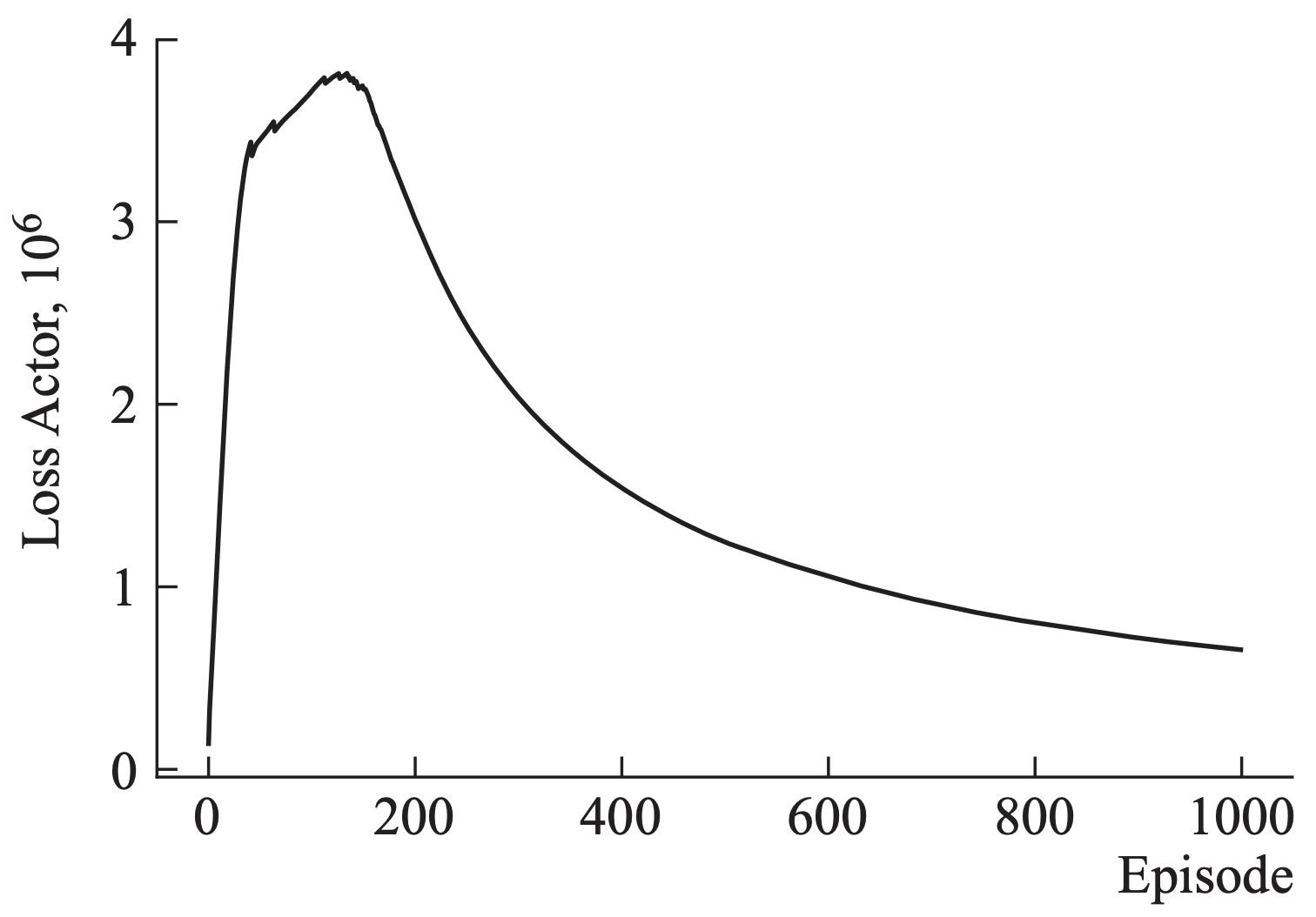}}
		\caption{The dependence of the loss function value on the episode of the Actor network.}
		\label{fig:actor_loss}
\end{figure}

    \begin{figure}[t]
\centering{\includegraphics[width=65mm]{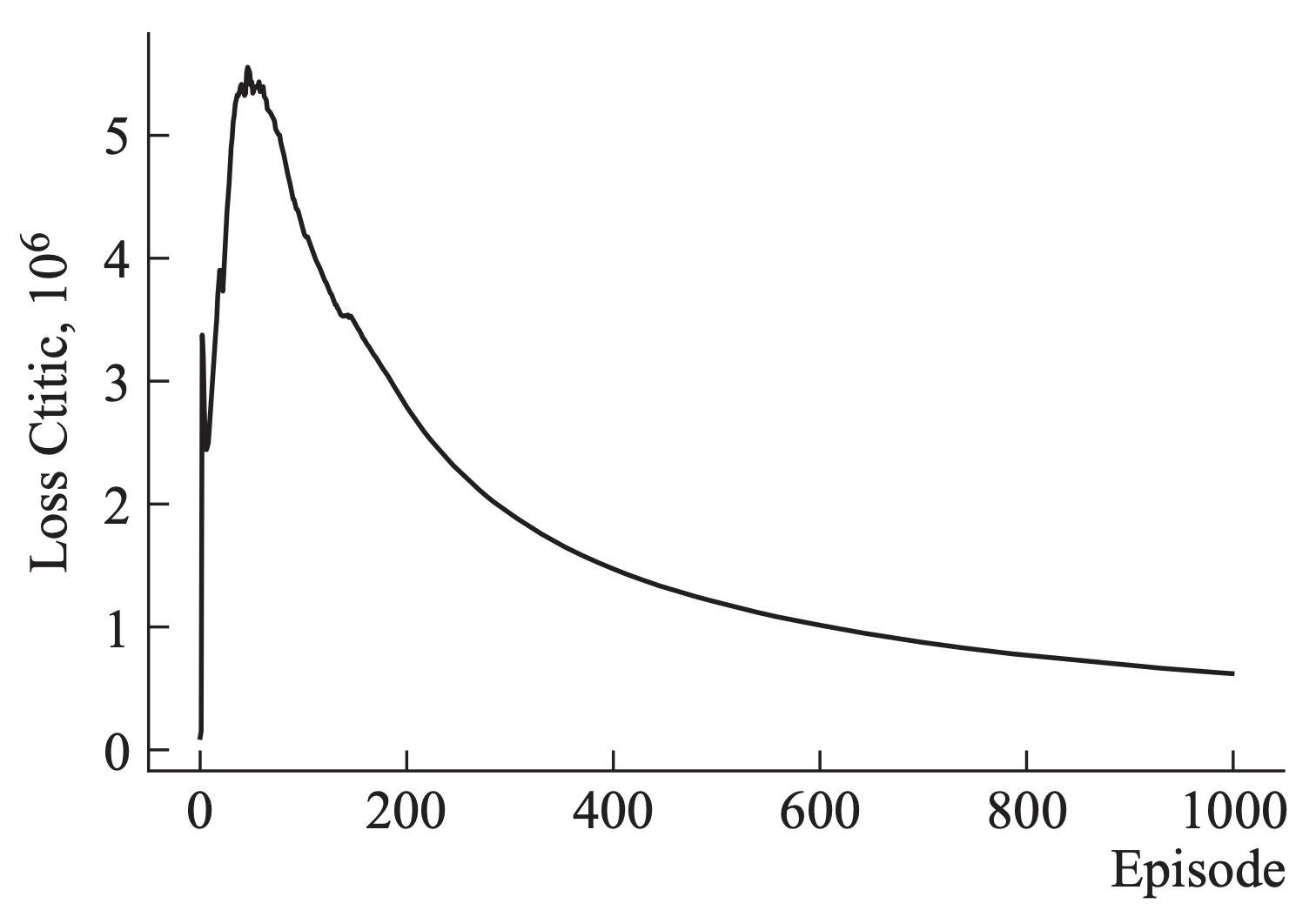}}
		\caption{The dependence of the loss function value on the episode of the Critic network.}
		\label{fig:critic_loss}
\end{figure}

Graphs of dependencies of the loss function of the Actor and Critic neural networks were also obtained. They are shown in Fig. \ref{fig:actor_loss} and \ref{fig:critic_loss} respectively.

    The graphs show a gradual decrease in the value of the loss function with an increase in training episodes, which indicates the correct choice of training coefficients.

\subsection{Learning result}

    \begin{figure}[t]
\centering{\includegraphics[width=130mm]{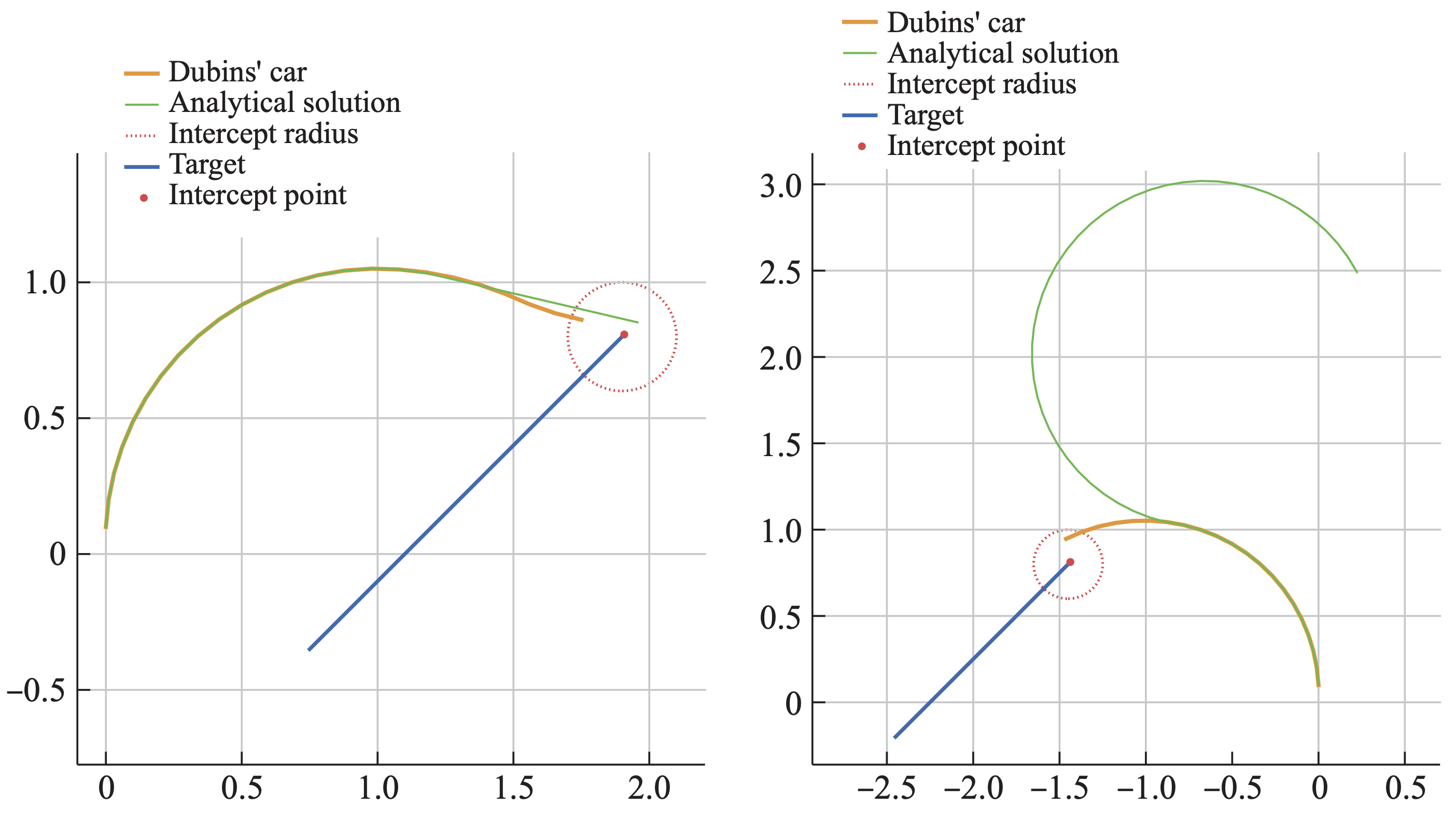}}
		\caption{Comparison of interception trajectories of a rectilinearly moving target with different initial parameters.}
		\label{fig:comp_tr_nn}
\end{figure}

In Fig. \ref{fig:comp_tr_nn} shows the trajectories obtained using a neural network and an analytical solution. The initial parameters of the target and the pursuer in this case had the values specified in Table \ref{table:4}.

\begin{table}[b]
	\caption{Initial parameters of the movement of the target and the pursuer\hspace*{11mm}}
    \label{table:4}
    \centering
    \begin{tabular}[H]{|p{8.2cm}|c|c|}
    \hline
    \multicolumn{1}{|c|}{\rule{0pt}{4mm}Parameter} & {Value} & {Value} \\ [0.25ex]
    \hline
    \rule{0pt}{4mm}{Initial coordinates of the target movement} & $(0{,}8; -0{,}4)$
     & $(-2{,}5; -0{,}25)$
    \\
     {Constant target rate $v_x$} & 0,5 & 0,5
    \\
     {Constant target rate $v_y$} & 0,5 & 0,5
    \\
     {Initial coordinates of the pursuer's movement} & $(0;0)$ & $(0;0)$
    \\
     {Initial orientation of the pursuer} $\varphi(0)$ & $\displaystyle\frac{\pi}{2}$ & $\displaystyle\frac{\pi}{2}$
    \\
    \small{Constant rate of the pursuer $v$} & 1 & 1
    \\
     {Intercept radius $\delta$} & 0,2 & 0,2
    \\
    \hline
    \end{tabular}
    \end{table}

long the trajectories shown in Fig. \ref{fig:comp_tr_nn}, it can be seen that the network was able to build
 a more efficient trajectory. In this case, the optimal interception time obtained using the analytical solution is
 ${T_{opt}\approx.5{,}42}$~s. And the time for which the network was able to intercept the target is
 ${T_{nn}\approx.2{,}1}$~s. This result is explained by the presence of the intercept radius $\delta=0{,}2$.

    \begin{figure}[t]
\centering{\includegraphics[width=130mm]{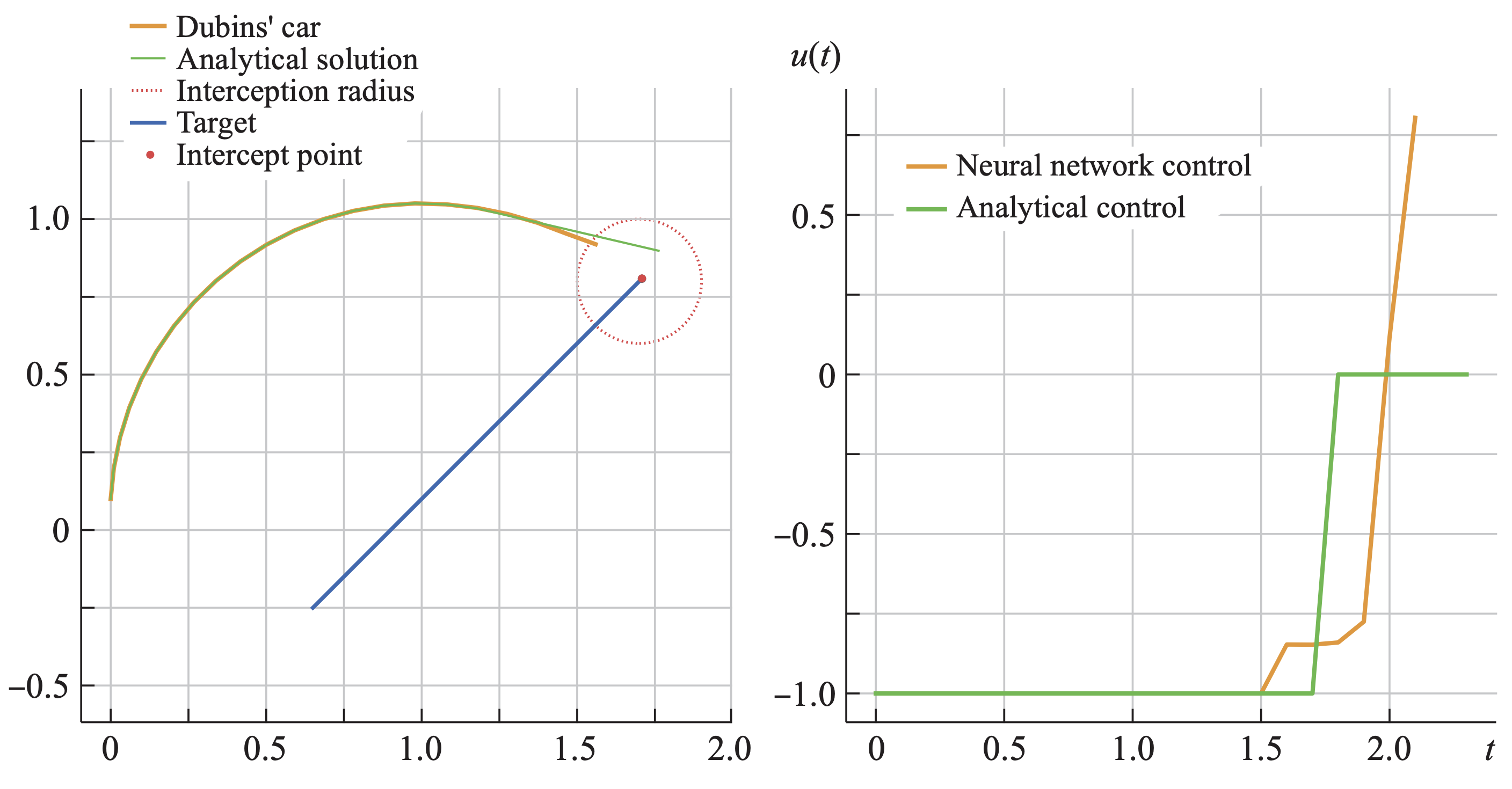}}
		\caption{Comparison of control functions from time.}
		\label{fig:comp_u}
\end{figure}

In Fig. \ref{fig:comp_u} on the right you can see a comparison of neural network control graphs with analytical. As can be seen, the controls differ significantly in the final section of the trajectory due to the fact that the neural network adjusts the terminal interception conditions. Optimal synthesis in a problem with an unfixed intercept angle consists of `Arc-line" or `Arc-arc" sections\cite{control}, and in a problem with a fixed intercept angle~--- in general, from the `Arc-line-arc" section \cite{Automatica}. It is the latter option that synthesizes the neural network. At the same time, as can be seen from Fig. \ref{fig:comp_tr_nn}, there is a section of the trajectory where the neural network chooses not the optimal, but close to the optimal value of the turning radius. The second reason for the difference is that the neural network optimizes the local reward function, which is different from the performance functional that was used when setting the task.

    \begin{figure}[t]
\centering{\includegraphics[width=130mm]{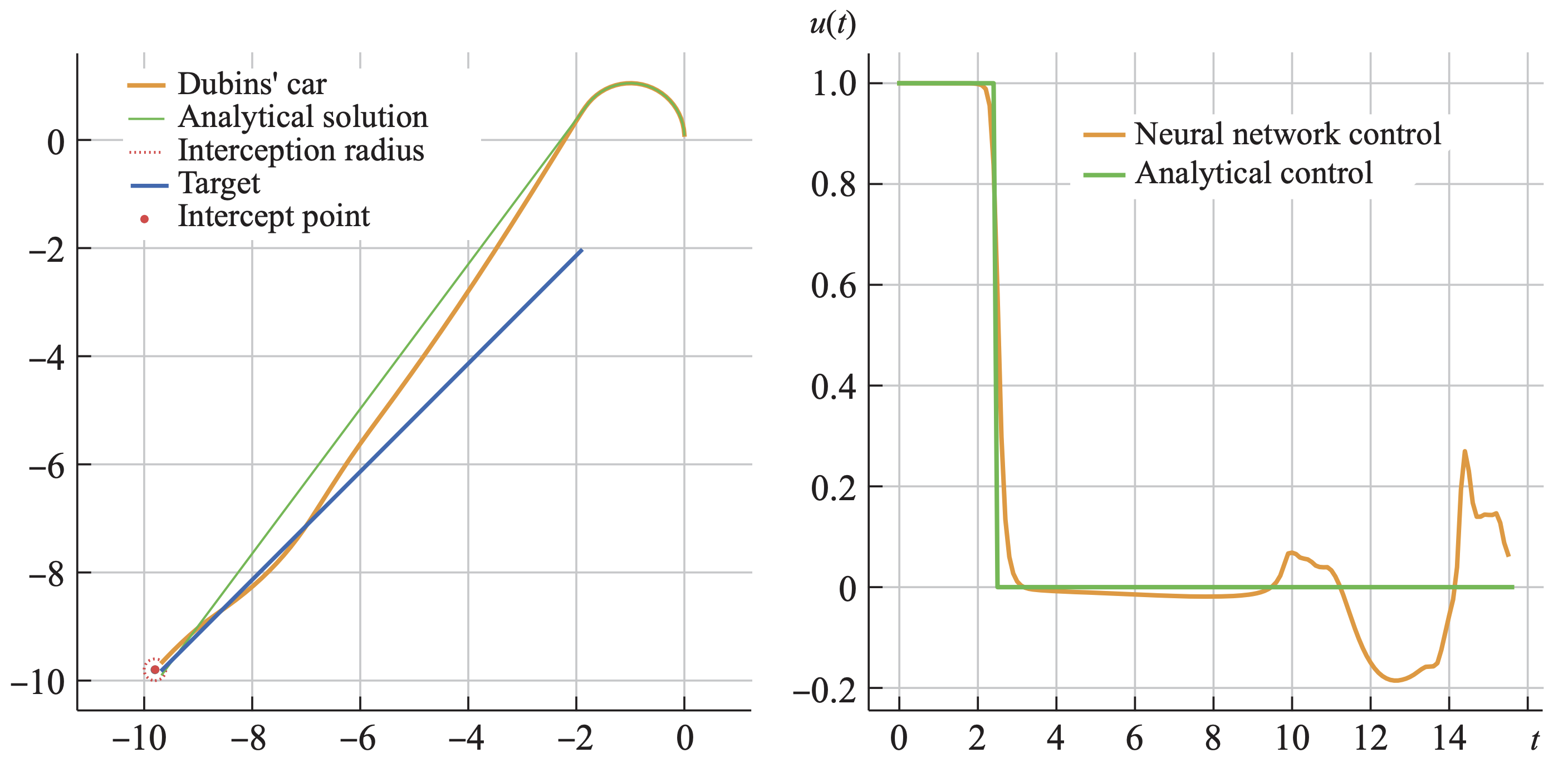}}
		\caption{Comparison of control functions from time.}
		\label{fig:comp_u_long}
\end{figure}

In Fig. \ref{fig:comp_u_long} shows the trajectories of intercepting the target and the dependence of the control function on time. On the right graph, it can be observed that the neural network control function has values close to optimal in the area where the analytical solution gives zero control. In addition, the deviations of the neural network control do not exceed the value of the intercept radius $\delta$. The interception times in this case are almost identical: $T_{opt}\approx T_{nn}\approx 21$ s.

    \subsection{Sensitivity analysis}
Let's analyze how much the neural network solution depends on the input parameters, since in theory the neural network should generalize the resulting solution well to states and parameters that it has not yet "seen" during training.

    \begin{figure}[t]
\centering{\includegraphics[width=130mm]{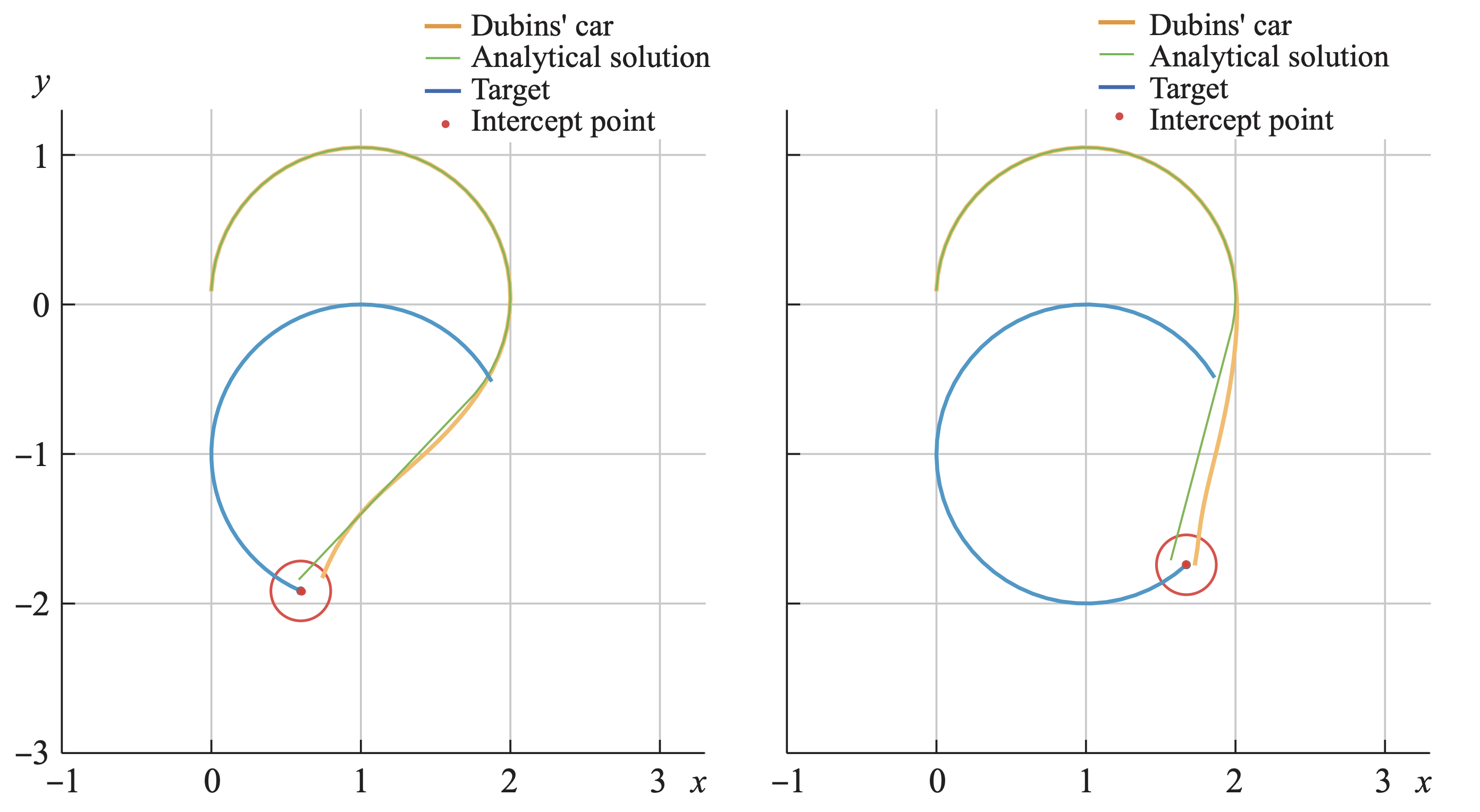}}
		\caption{Comparison of circular intercept trajectories at different initial parameters.}
		\label{fig:compofw}
\end{figure}

    \begin{figure}[t]
\centering{\includegraphics[width=130mm]{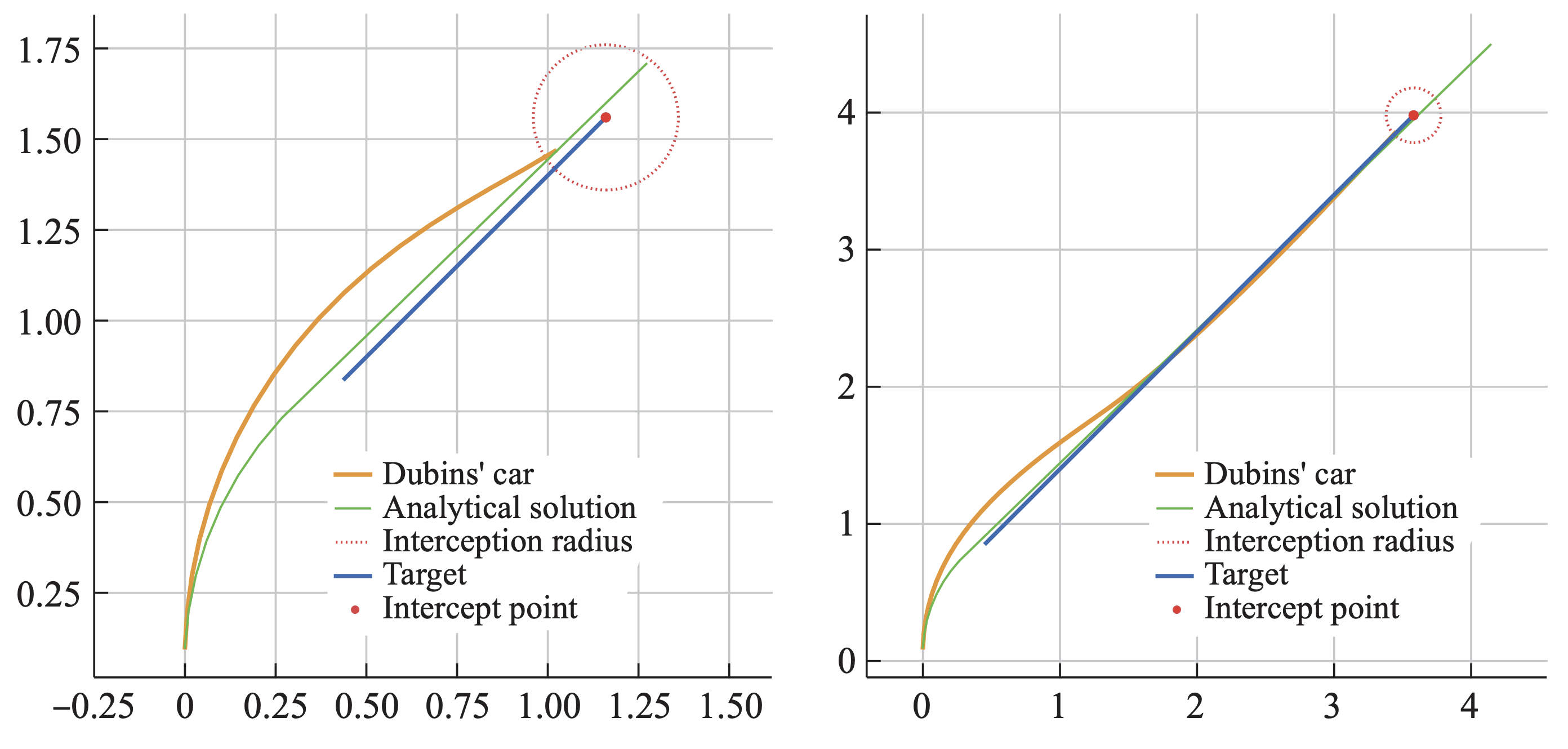}}
		\caption{Comparison of target intercept trajectories at different initial parameters.}
		\label{fig:comp_tr_velo}
\end{figure}

To intercept a target moving in a circle, we will train the neural network only on targets with a single radius and a single angular velocity and check whether it can successfully catch a target with other parameters. As can be seen in Fig. \ref{fig:compofw}, the network successfully copes with the task, in the left figure the angular velocity of the target is 0,7 of the angular velocity used in training, in the right figure the interception of an ordinary target is depicted. Experiments were conducted for angular velocity values from 0,7 to 1,3, in which the neural network successfully intercepted the target.

In the case of interception of a rectilinearly moving target, the neural network was trained at the target velocity values ${v_x = v_y = 0{,}5}$. In Fig. \ref{fig:comp_tr_velo} shows the results of network testing with speeds differing by $20\%$~--- in the left figure, the target has a speed of ${v_x = v_y =0{,}4}$, and on the right ${v_x = v_y =0{,}6}$.

        It follows from the results obtained that the network generalizes the solution well. This can be useful for applied tasks, since in them the parameters are often known with some error.

\section{Conclusion}

The paper proposed two DDPG-based neural network algorithms for the synthesis of trajectories of interception by the Dubins' car of targets moving along rectilinear and circular trajectories. The features of the proposed algorithms are their ability to work with the space of continuous actions, the guarantee of learning and working with different relative initial positions of goals and the Dubins' car. It is shown that the network successfully generalizes the solution and in some situations offers the fastest solution to the interception problem.

The undoubted advantages of the proposed algorithms can be used, and the algorithms themselves are modified to obtain a barrier surface in the differential game of two cars.



\bibliography{sn-bibliography.bib}


%


\end{document}